\documentclass[a4paper]{article}
\usepackage{amsmath}
\usepackage{amsbsy}
\usepackage{amsfonts}
\usepackage{amstext}
\usepackage{amssymb}
\usepackage{mathbbol}
\usepackage{wasysym}
\usepackage{eucal} 
\usepackage{mathtools}
\usepackage[amsmath,thmmarks]{ntheorem}		
\usepackage{stackrel}
\usepackage{youngtab}

\usepackage{tikz}
   \usetikzlibrary{matrix}
   \usetikzlibrary{calc}
   \usetikzlibrary{decorations,snakes}

\usepackage{float}
\usepackage{subfigure}

\usepackage{enumitem}
\usepackage{mathtools}

\usepackage[a4paper, left=2.5cm, right=2.5cm, top=3.0cm, bottom=2.5cm, bindingoffset=1cm]{geometry}
\usepackage{color}

\newcommand{\ground}{\ensuremath{\mathbb k}} 

\newcommand{\ZZ}{\ensuremath{{\mathbb Z}}}
\newcommand{\CC}{\ensuremath{{\mathbb C}}}

\newcommand{\Mod}{\operatorname{mod}}

\newcommand{\End}{\ensuremath{\text{End}}}

\renewcommand{\dim}{\ensuremath{\text{dim}}}

\renewcommand{\i}{\ensuremath{^{-1}}}

\newcommand{\p}{\ensuremath{^{\prime}}}
\newcommand{\pp}{\ensuremath{^{\prime\prime}}}
\newcommand{\ppp}{\ensuremath{^{\prime\prime\prime}}}

\renewcommand{\t}[1]{\ensuremath{\widetilde{#1}}}

\newcommand{\h}[1]{\ensuremath{\widehat{#1}}}

\newcommand{\ul}[1]{\ensuremath{\underline{#1}}}

\newcommand{\rar}{\rightarrow}

\newcommand{\ty}[1]{\ensuremath{\mathsf{#1}}} 											

\newcommand{\liesl}{\ensuremath{{\mathfrak s \mathfrak l}}}

\newcommand{\qq}{\ensuremath{q}}

\newcommand{\basi}{\ensuremath{{\cal B}}}														

\newcommand{\plac}[1]{\ensuremath{{\cal P}_{#1}}} 									
\newcommand{\partic}[1]{\ensuremath{{\cal P}^\text{part}_{#1}}} 		
\newcommand{\aplac}[1]{\ensuremath{\h{\cal P}_{#1}}} 								

\newcommand{\bosconfi}[1]{\ensuremath{\ul{\mathbf{\textsc{#1}}}}}

\theoremstyle{plain}
\theoremseparator  {.}           
\theoremheaderfont {\bfseries}
\theorembodyfont   {\normalfont} 
\theoremnumbering  {arabic}
\makeatletter
\newtheoremstyle{normal}%
{\item[\hskip\labelsep \theorem@headerfont ##1\ ##2\theorem@separator]\normalfont}%
{\item[\hskip\labelsep \theorem@headerfont ##1\ ##2]{\theorem@headerfont (##3)}\theorem@separator\ \normalfont}
\newtheoremstyle{nonumber}%
{\item[\theorem@headerfont\hskip\labelsep ##1\theorem@separator]\normalfont}%
{\item[\theorem@headerfont\hskip \labelsep ##1]{\theorem@headerfont (##3)}\theorem@separator\ \normalfont}
\makeatother

\theoremstyle{normal}
\newtheorem{thm}               {Theorem} [section]
\newtheorem{theo}       [thm]  {Theorem}              
\newtheorem{lemma}      [thm]  {Lemma}             
\newtheorem{cor}        [thm]  {Corollary}

\newtheorem{defi}       [thm]  {Definition}

\newtheorem{prop}       [thm]  {Proposition}
\newtheorem{bem}        [thm]  {Remark}
\newtheorem{bsp}        [thm]  {Example}

\theoremstyle{nonumber}
\theoremsymbol{$\square$}     
\newtheorem{bew}{Proof}
\theoremsymbol{}              

\theoremstyle{nonumber}
\newtheorem{thmnn}{Theorem}

\begin{document}

\title{A plactic algebra action on bosonic particle configurations: The classical case}
\author{Joanna Meinel}
\date{}
\maketitle
\begin{abstract}
We study the plactic algebra and its action on bosonic particle configurations in the classical case. These particle configurations together with the action of the plactic generators can be identified with crystals of the quantum analogue of the symmetric tensor representations in type $\ty{A}$.  It turns out that this action factors over a quotient algebra that we call partic algebra, whose induced action on bosonic particle configurations is faithful. We describe a basis of the partic algebra explicitly in terms of a normal form for monomials, and we compute the center of the partic algebra.
\end{abstract}

\section*{Introduction}

\setcounter{section}{0}

Let $\ground$ be a field. Fix an integer $N\geq3$.
The (local) plactic algebra $\plac{N}$ is the unital associative $\ground$-algebra generated by $a_1,\ldots,a_{N-1}$ subject to the plactic relations
\begin{align}
\label{plac1} a_i a_{i-1} a_i\ &=\ a_i a_i a_{i-1}\quad &\text{for } 2\leq i\leq N-1,\\
\label{plac2} a_i a_{i+1} a_i\ &=\ a_{i+1} a_i a_i\quad &\text{for } 1\leq i\leq N-2,
\end{align}
together with the commutativity relation
\begin{align}
\label{comm} a_i a_j\ &=\ a_j a_i\quad &\text{for }|i-j|>1.
\end{align}
Our results hold over an arbitrary unitary associative ring, if we adapt the notation, replacing vector spaces by free modules and so on. In statements about the center we need to assume commutativity of the ground ring. For simplicity we choose to work over a field.

The plactic relations go back to Lascoux and Sch{\"u}tzenberger \cite{ls-plac}. They study the monoid defined by the ``plaxic relations'' (in the original, ``plaxique'' or ``a placche'') \eqref{plac1}, \eqref{plac2} and the non-local Knuth relation, a slightly weaker commutativity relation  $(a_i a_j) a_k = (a_j a_i) a_k$, $a_k (a_i a_j) = a_k (a_j a_i)$ for $i<k<j$ (in particular for $|i-j|>1$). This monoid is isomorphic to the monoid of semistandard Young tableaux with entries $1,\ldots,N-1$ (and multiplication defined by row bumping) by reading off the entries of a tableau from left to right and bottom to top, see \cite[Section~2.1]{fulton-yt} for the details. 
\par
The name local plactic algebra for the algebra defined by the relations~\eqref{plac1}, \eqref{plac2} and \eqref{comm} goes back to \cite{fomingreene} due to the additional ``local'' commutativity relation~\eqref{comm}. Fomin and Greene develop a theory of Schur functions in noncommutative variables that applies in particular to the (local) plactic algebra, see \cite[Example 2.6]{fomingreene}, including a generalized Littlewood-Richardson rule for Schur functions defined over the plactic algebra. The plactic algebra acts on Young diagrams by Schur operators, i.e. $a_i$ adds a box in the $i$th column if possible, and otherwise maps the diagram to zero \cite{fomin-schur}.
\par
The monoid defined by the plactic relations~\eqref{plac1}, \eqref{plac2} and \eqref{comm} appears as a Hall monoid or ``quantic monoid'' of type $\ty{A}_{N-1}$ in  \cite{reinge}, \cite{reinmono}:
Reineke defines the structure of a monoid on isomorphism classes of modules over the path algebra $\ground Q$ for an oriented Dynkin quiver $Q$. The product of two isomorphism classes $[M]$ and $[M\p]$ is defined by $[M\ast M\p]$, the isomorphism class of the generic extension of $M$ by $M\p$ in $\ground Q-\Mod$. The generic extension is up to isomorphism uniquely determined to be the extension with $\dim(\End_{\ground Q}(M\ast M\p))$ minimal among all possible extensions. Equivalently, the orbit of the generic extension is dense in the subset of extensions of $[M]$ and $[M\p]$ inside the representation variety of $Q$.
In particular, Reineke shows that for $Q=\ty{A}_{N-1}$ (with orientation given e.g. by $i\rar (i-1)$ for the vertices $2\leq i\leq N-1$ of $Q$), the $\ground$-linearisation of the resulting monoid is isomorphic the plactic algebra as defined above, where the isomorphism classes of the one-dimensional simple modules $[S_i]$ are mapped to the generators $a_i$. This is furthermore identified with the positive half of the twisted quantum group at $\qq=0$, which is obtained by twisting the multiplication and desymmetrizing the quantum Serre relations so that they can be rewritten without appearance of $\qq\i$.
In \cite{reinmono} it is proven that the twisted (positive) half of the quantum group specialized to $\qq=0$ is isomorphic to the linearisation of the Hall monoid. The desymmetrized Serre relations at $\qq=0$ are the plactic relations.
By Ringel's theorem \cite{ringel} we know that the positive half of the twisted quantum group is isomorphic to the generic Hall algebra for any Dynkin quiver $Q$. Hence, the specialisation of the generic Hall algebra at $\qq=0$ gives the Hall monoid. Different normal forms for monomials in the plactic algebra are given in terms of enumerations of the roots \cite[Theorem~2.10]{reinmono}.
\par
In \cite{ks} the plactic algebra appears in the study of bosonic particle configurations. Schur functions in the generators of the affine plactic algebra are defined using Bethe Ansatz techniques to show that they are well-defined despite the noncommutativity of the generators. Combinatorially, a bosonic particle configuration is given by a tuple $(k_1,\ldots,k_{N-1},k_0)$ in $\ZZ_{\geq0}^N$. One can think of such a tuple as a finite number of particles distributed on a discrete lattice of $N$ positions on a line segment (the classical case) or along a circle (the affine case). Here we focus on the classical case. The generators $a_i$ act on the particle configurations and their $\ground$-span by lowering $k_i$ by $1$ and increasing $k_{i+1}$ by $1$, if possible. If not possible since $k_i=0$, the result is $0$. In the picture this would correspond to propagation of a particle from position $i$ to $i+1$.
One can identify bosonic particle configurations with Young diagrams, then the operator $a_i$ acts by adding a box in the $(i+1)$-st row of the Young diagram, if possible, and by $0$ otherwise. Up to an index shift and switching rows and columns, this is the same as the action on Young diagrams by Schur operators from \cite{fomingreene}. We will use the identification of the $\ground$-span of bosonic particle configurations with the vector space of polynomials $\ground[x_1,\ldots,x_{N-1},x_0]$ so that a particle configuration $(k_1,\ldots,k_{N-1},k_0)$ corresponds to a monomial $x_1^{k_1}\ldots x_{N-1}^{k_{N-1}}x_0^{k_0}$. Then the generator $a_i$ of the plactic algebra acts by lowering the exponent of $x_i$ by $1$ and raising the exponent of $x_{i+1}$ by $1$. Note that this action is combinatorial in the sense of \cite{fomingreene}.
\begin{figure}[H]
\begin{center}
\begin{tikzpicture}[scale=0.6]
\begin{scope}[xshift=0cm]
\path[use as bounding box] (-4.5,-1.5) rectangle (4.5,2.5);
\draw[thick] (-4.3,0) -- (3.3,0);

\foreach \x in {1,...,8}
\draw[thick] ({-5+\x},-0.1) -- ({-5+\x},0.1);


\foreach \x in {1,...,7}
\draw ({-5+\x},-0.4) node {$\x$};

\draw (3,-0.4) node {$0$};

\foreach \x in {1,4,6,7,8}
\fill ({-5+\x},0.5) circle (1.2mm);
\foreach \x in {1,7,8}
\fill ({-5+\x},0.9) circle (1.2mm);
\foreach \x in {1}
\fill ({-5+\x},1.3) circle (1.2mm);
\foreach \x in {1}
\fill ({-5+\x},1.7) circle (1.2mm);
\end{scope}
\end{tikzpicture}
\end{center}
\caption{Example for $N=8$: A bosonic particle configuration on a line segment given by the tuple $(4,0,0,1,0,1,2,2)$.}
\label{fig:intro}
\end{figure}
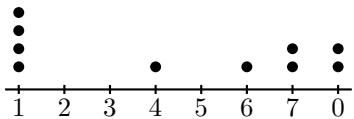
A very prominent combinatorial realization of the action of the generators $a_i$ of the plactic algebra on bosonic particle configurations with $k$ particles is the action of the Kashiwara operators $\t{f}_i$ on the crystal ${\cal{B}}(k\omega_1)$ of type $\ty{A}_{N-1}$, i.e. the crystal for the quantum analogue of the symmetric representation $L(k\omega_1)=\text{Sym}^k(\CC^N)$ of the Lie algebra $\liesl_N(\CC)$ (see e.g. \cite{hongkang} for details). The Young tableaux that constitute the nodes in the crystal graph correspond to the particle configurations by placing a particle at position $i$ for each box labelled $i$ in the Young tableau. Similarly, the crystal $\basi(\omega_k)$ for the quantum analogue of the alternating representation can be identified with fermionic particle configurations together with the action of the generators of the plactic algebra. This implies that the Kashiwara operators satisfy the plactic relations on $\basi(\omega_k)$ and $\basi(k\omega_1)$, see also \cite[Chapter I.1]{diss}. Relations among Kashiwara operators for abstract crystals of simply laced finite and affine type were studied by Stembridge in \cite{stembridge} where a list of relations is given that hold if and only if the abstract crystal graph can be realized as a crystal graph of an integrable highest weight representation.
\par
Here we study the representation of the plactic algebra on bosonic particle configurations more closely. Our main goal is to identify the kernel of this representation and to describe the resulting algebra. For Young diagrams our results can be interpreted as a full list of generating relations among the Schur operators. For crystals they can be interpreted as relations satisfied by the Kashiwara operators $\t{f}_i$ on crystals of the form ${\cal{B}}(k\omega_1)$.
\par
Let us point out that from \cite[Proposition 9.1]{ks}, \cite[Proposition 2.4.1]{bfz}, \cite{bjs} it is known that on fermionic particle configurations the nilTemperley-Lieb quotient of the plactic algebra acts faithfully. In \cite{bejo} the case of affine fermionic particle configurations was studied, including a description of a normal form for monomials in the affine nilTemperley-Lieb algebra and its center.
\par \ 
\par
\textbf{Organization and Results.} 
In Section~\ref{sec:classipartic} we introduce a quotient of the classical plactic algebra named partic algebra, and we start by some small technical preliminaries.
In Section~\ref{sec:particbasis} we construct a normal form of the monomials in the partic algebra, leading to the Basis Theorem \ref{theo:particbasis}:
\begin{thmnn}
The partic algebra $\partic{N}$ has a basis given by monomials of the form
\begin{align*}
\{ a_{N-1}^{d_{N-1}}\ldots a_2^{d_2} a_1^{k_1}a_2^{k_2}\ldots a_{N-1}^{k_{N-1}}\ |\ d_i\leq d_{i-1}+k_{i-1}\ \text{ for all }\ 3\leq i\leq N-1,\ d_2\leq k_1\}
\end{align*}
where $d_i,k_i\in\ZZ_{\geq0}$ for all $1\leq i\leq N-1$.
\end{thmnn}
In Section~\ref{sec:bosonicac} we discuss the action of the classical plactic and the partic algebra on bosonic particle configurations which we realize as an action on the polynomial ring $\ground[x_1,\ldots,x_{N-1},x_0]$, and we prove faithfulness of the action of the partic algebra in Theorem~\ref{theo:particfaith}:
\begin{thmnn}
The action of the partic algebra $\partic{N}$ on $\ground[x_1,\ldots,x_{N-1},x_0]$ is faithful.
\end{thmnn}
In Section~\ref{sec:bosce} we describe the center of the partic algebra in terms of Theorem \ref{theo:centpart}:
\begin{thmnn}
The center of the partic algebra $\partic{N}$ is given by the $\ground$-span of the elements
$$\{ a_{N-1}^r a_{N-2}^r\ldots a_2^r a_1^r\ |\ r\geq0\}.$$
\end{thmnn}
In Section~\ref{sec:aff} we give an outlook to the much harder affine case treated in the followup work \cite{affplac}. We recall the definition of the affine partic algebra and its action on affine bosonic particle configurations. For the description of the kernel we find an unexpected generalization of the partic relation from the classical case.
\par \ 
\par
This paper goes back to the author's PhD thesis: Most of Sections~\ref{sec:classipartic},~\ref{sec:particbasis},~\ref{sec:bosonicac} and~\ref{sec:bosce} can be found in \cite[Chapter I.3]{diss}, while Section~\ref{sec:aff} is new and replaces the outdated view from \cite[Chapter I.3.6]{diss}. 
\par \ 
\par
\textbf{Acknowledgements.} I would like to thank the Max Planck Institute for Mathematics and the Hausdorff Center for Mathematics for excellent research conditions and funding the project through IMPRS/BIGS. I thank Catharina Stroppel for supervising the thesis, and I am grateful to Michael Ehrig and Daniel Tubbenhauer for helpful discussions.

\section{The partic algebra}\label{sec:classipartic}
We introduce a quotient of the classical plactic algebra:
\begin{defi}
Define the partic algebra $\partic{N}$ to be the quotient of $\plac{N}$ by the additional relations
\begin{align}
\label{part}
a_i a_{i-1} a_{i+1} a_i \ &=\ a_{i+1} a_i a_{i-1} a_i\quad &\text{for } 2\leq i\leq N-2.
\end{align}
\end{defi}
Note that one can interpret the plactic relations~\eqref{plac1}, \eqref{plac2} as commutativity of the product $(a_{i+1}a_i)$ with the generators $a_{i+1}$ and $a_i$.
Relation~\eqref{part} together with~\eqref{plac1} implies in particular that $(a_{i+1}a_i)$ and $(a_ia_{i-1})$ commute.
\begin{bem}
This relation appears naturally in the study of bosonic particle configurations, see Section~\ref{sec:bosonicac}.
In contrast, in the Hall monoid of finite type $\ty{A}_{N-1}$ one cannot expect $[S_{i+1}\ast S_i]$ and $[S_i\ast S_{i-1}]$ to commute. This is because precisely one of $S_{i+1}\ast S_i\ast S_i\ast S_{i-1}$, $S_i\ast S_{i-1}\ast S_{i+1}\ast S_i$ is a nontrivial extension of $S_{i+1}\ast S_i$ and $S_i\ast S_{i-1}$ (it depends on the choice of orientation which one is nontrivial) -- in much the same way as $[S_i]$ and $[S_{i\pm1}]$ do not commute.
\end{bem}
\begin{bem}\label{bem:particgrading}
We have two gradings on both the plactic and the partic algebra:
\begin{enumerate}
\item
All relations preserve the length of monomials, hence $\plac{N}$ and $\partic{N}$ can be equipped with a $\ZZ$-grading by the length of monomials.
\item
All relations preserve the number of different generators in a monomial, hence $\plac{N}$ and $\partic{N}$ can be equipped with a $\ZZ^{N-1}$-grading that assigns to the generator $a_i$ the degree $e_i$, the $i$-th standard basis vector in $\ZZ^{N-1}$. This is a refinement of the above length grading.
\end{enumerate}
\end{bem}
\begin{lemma}\label{lem:plactitech} In the plactic (and hence also in the partic) algebra, the following relations hold:
\begin{enumerate}[label=\roman{*}), ref=\ref{lem:plactitech}.(\roman{*})]
\item\label{lem:plactitech1}For all generators $a_i$, $a_{i-1}$, $2\leq i\leq N-1$ and all $m\geq0$, we have
\begin{align}
a_i^m a_{i-1}^m\ &=\ (a_i a_{i-1})^m,\label{power}\\
a_i(a_i^m a_{i-1}^m)\ &=\ (a_i^m a_{i-1}^m) a_i.\notag
\end{align}
\item\label{lem:plactitech2}
For all $i\geq k\geq j$ we have
\begin{align}\label{wurmcomm}
(a_i a_{i-1}\ldots a_{j+1} a_j)a_k\ =\ a_k(a_ia_{i-1}\ldots a_{j+1}a_j).
\end{align}
\end{enumerate}
\end{lemma}
\begin{bew}
\begin{enumerate}
\item
The second equation of Lemma~\ref{lem:plactitech1} follows from the first by the plactic relation~\eqref{plac1}.
By induction, $a_i^m a_{i-1}^m = a_i (a_i a_{i-1})^{m-1} a_{i-1}= (a_i a_{i-1})^{m-1}a_i a_{i-1} = (a_i a_{i-1})^m$.
\item
This equality follows from the calculation
\begin{align*}
(a_i a_{i-1}\ldots a_{j+1} a_j)a_k\ &\stackrel{\eqref{comm}}{=}\ a_ia_{i-1} \ldots a_{k+1}a_k a_{k-1}a_k\ldots a_{j+1}a_j\\
&\stackrel{\eqref{plac1}}{=}\ a_ia_{i-1} \ldots a_{k+1}a_k a_k a_{k-1}\ldots a_{j+1}a_j\\
&\stackrel{\eqref{plac2}}{=}\  a_ia_{i-1} \ldots a_k a_{k+1} a_k a_{k-1}\ldots a_{j+1}a_j\\
&\stackrel{\eqref{comm}}{=}\ a_k(a_ia_{i-1}\ldots a_{j+1}a_j).
\end{align*}
\end{enumerate}
\end{bew}

\section{A basis of the partic algebra}\label{sec:particbasis}
In this section we formulate the following main theorem:
\begin{theo}\label{theo:particbasis}
The partic algebra $\partic{N}$ has a basis given by monomials of the form
\begin{align}
\label{cond}\{ a_{N-1}^{d_{N-1}}\ldots a_2^{d_2} a_1^{k_1}a_2^{k_2}\ldots a_{N-1}^{k_{N-1}}\ |\ d_i\leq d_{i-1}+k_{i-1}\ \text{ for all }\ 3\leq i\leq N-1,\ d_2\leq k_1\}
\end{align}
where $d_i,k_i\in\ZZ_{\geq0}$ for all $1\leq i\leq N-1$.
\end{theo}
Our approach is based on the observation that it suffices to construct a normal form for monomials to obtain a $\ground$-basis for an algebra given by generators and monomial relations. This follows from \cite[Proposition 2]{hey} where it is proven that an algebra of the form 
$$\ground[X]/( m-m\p\ |\ m,m\p\text{ certain words in the generators in }X)$$
is isomorphic to the semigroup algebra $\ground[S]$ for the semigroup $S$ defined by the same set of generators $X$ and relations $m=m\p$ for all $m-m\p$ that generate the ideal. By definition $\ground[S]$ has a basis given by elements in $S$ which can be represented by monomials in $X$.
\begin{bem}
This approach can be extended to algebras of the form 
$$A=\ground[X]/( m,m\p-m\pp\ |\ m,m\p,m\pp\text{ certain words in the generators in }X).$$ 
The corresponding semigroup $S$ is given by generators $X\cup\{\emptyset\}$, where $\emptyset$ denotes the absorbing element, and the relations are of the form $m=\emptyset$, $m\p=m\pp$ where $m$, $m\p-m\pp$ are generators of the ideal. Then $A\cong\ground[S]/(\emptyset)$, and a basis of $A$ is given by a normal form for the elements in $S\setminus\{\emptyset\}$ (see \cite[Remark I.2.5.1]{diss}).
\end{bem}
In this section we show that every monomial in the partic algebra is equivalent to a monomial of the form~\eqref{cond}. In Section~\ref{sec:bosonicac} we observe that these monomials act pairwise differently on the particle configuration module, and we conclude that they must have been distinct.
\begin{prop}\label{prop:nofo}
Every monomial in the partic algebra $\partic{N}$ is equivalent to a monomial of the form~\eqref{cond}, i.e. $a_{N-1}^{d_{N-1}}\ldots a_2^{d_2} a_1^{k_1}a_2^{k_2}\ldots a_{N-1}^{k_{N-1}}$ with $d_i\leq d_{i-1}+k_{i-1}$ for all $3\leq i\leq {N-1}$ and $d_2\leq k_1$.
\end{prop}
\begin{bew}
The proof works by induction on the length of monomials.
If the length is equal to $1$, we have $a_i=a_i^{k_i}$ for $k_i=1$, and the condition from~\eqref{cond} is preserved. For the induction step our goal is to show that
\begin{align}\label{goal}
a_i\cdot\left(a_{N-1}^{d_{N-1}}\ldots a_2^{d_2} a_1^{k_1}a_2^{k_2}\ldots a_{N-1}^{k_{N-1}}\right)\ =\ a_{N-1}^{d_{N-1}}\ldots a_i^{d_i\p}\ldots a_2^{d_2} a_1^{k_1}a_2^{k_2}\ldots a_i^{k_i\p}\ldots a_{N-1}^{k_{N-1}},
\end{align}
where $d_i\p=d_i$ and $k_i\p=k_i+1$, or $d_i\p=d_i+1$ and $k_i\p=k_i$ are such that the inequality condition~\eqref{cond} is preserved.
Since we can commute $a_i$ with all $a_j$ as long as $j\neq i\pm1$, we only need to consider 
$$a_i\cdot\left(a_{i+1}^{d_{i+1}}a_i^{d_i} a_{i-1}^{d_{i-1}}\ldots a_{i-1}^{k_{i-1}}a_i^{k_i}\right).$$
In order to prove that this can be rewritten as in~\eqref{goal}, we have to show that either we can pass $a_i$ through to the right hand side, increasing the exponent $k_i$ by one, or we leave it at the left hand side, increasing $d_i$ by one.
\begin{enumerate}
\item Case $d_{i+1}=d_i=d_{i-1}=k_{i-1}=0$:
Set $k_i\p=k_i+1$. The inequality condition~\eqref{cond} is automatically satisfied if we increase one of the $k$'s, so there is nothing to check. The equality ~\eqref{goal} is obvious since we only apply the commutativity relation~\eqref{comm}.
\item Case $d_{i+1}=d_i=d_{i-1}=0$, $k_{i-1}>0$:
Set $d_i\p=1$. The inequality condition~\eqref{cond} is preserved since $k_{i-1}\geq1$, and again we only apply the commutativity relation~\eqref{comm}.
\item Case $d_{i+1}=d_i=0$,  $d_{i-1}>0$, $k_{i-1}$ arbitrary:
Set $d_i\p=1$. The inequality condition~\eqref{cond} is preserved since $d_{i-1}\geq1$, and as before we only apply the commutativity relation~\eqref{comm}.
\item Case $d_{i+1}=0$, $d_i>0$, $d_{i-1}$ and $k_{i-1}$ arbitrary so that $d_i\leq d_{i-1}+k_{i-1}$:
\begin{itemize}
\item $d_i< d_{i-1}+k_{i-1}$: Set $d_i\p=d_i+1$.
\item $d_i= d_{i-1}+k_{i-1}$: We cannot increase $d_i$, hence we have to show that we can commute $a_i$ past $a_{i-1}^{d_{i-1}}$ and $a_{i-1}^{k_{i-1}}$ to increase $k_i$. Indeed, we can apply equality~\eqref{power} from Lemma~\ref{lem:plactitech1} to obtain
\begin{align*}
a_i\left(a_i^{d_i} a_{i-1}^{d_{i-1}}\ldots a_{i-1}^{k_{i-1}}a_i^{k_i}\right)\ &=\ a_i^{d_{i-1}+k_{i-1}+1} a_{i-1}^{d_{i-1}}\ldots a_{i-1}^{k_{i-1}}a_i^{k_i}\\
&=\ a_i^{k_{i-1}+1} (a_i a_{i-1})^{d_{i-1}}\ldots a_{i-1}^{k_{i-1}}a_i^{k_i}\\
&=\ (a_i a_{i-1})^{d_{i-1}}a_i^{k_{i-1}+1}\ldots a_{i-1}^{k_{i-1}}a_i^{k_i}\\
&=\ (a_i a_{i-1})^{d_{i-1}}\ldots a_i^{k_{i-1}+1} a_{i-1}^{k_{i-1}}a_i^{k_i}\\
&=\ (a_i a_{i-1})^{d_{i-1}}\ldots a_i(a_i a_{i-1})^{k_{i-1}}a_i^{k_i}\\
&=\ (a_i a_{i-1})^{d_{i-1}}\ldots (a_i a_{i-1})^{k_{i-1}}a_i^{k_i+1}\\
&=\ (a_i a_{i-1})^{d_{i-1}}a_i^{k_{i-1}}\ldots a_{i-1}^{k_{i-1}}a_i^{k_i+1}\\
&=\ a_i^{d_{i-1}+k_{i-1}} a_{i-1}^{d_{i-1}}\ldots a_{i-1}^{k_{i-1}}a_i^{k_i+1}.
\end{align*}
\end{itemize}
\item Case $d_{i+1}>0$, $d_{i}$, $k_{i}$ and $d_{i-1}$, $k_{i-1}$ arbitrary so that $d_{i+1}\leq d_{i}+k_{i}$, $d_i\leq d_{i-1}+k_{i-1}$:
We reduce to the previous cases by proving
$$a_ia_{i+1}^{d_{i+1}}\cdot\left(a_i^{d_i} a_{i-1}^{d_{i-1}}\ldots a_{i-1}^{k_{i-1}}a_i^{k_i}\right)\ =\ a_{i+1}^{d_{i+1}}a_i\left(a_i^{d_i} a_{i-1}^{d_{i-1}}\ldots a_{i-1}^{k_{i-1}}a_i^{k_i}\right).$$
\begin{itemize}
\item $d_{i+1}\leq d_i$:
Here we can apply Lemma~\ref{lem:plactitech1} to obtain
$$a_i(a_{i+1}^{d_{i+1}}a_i^{d_i})\ =\ a_i (a_{i+1}a_i)^{d_{i+1}}a_i^{d_i-d_{i+1}}\ =\ (a_{i+1}a_i)^{d_{i+1}}a_i^{d_i-d_{i+1}}\ =\ a_{i+1}^{d_{i+1}}a_i a_i^{d_i}.$$
\item $d_{i+1}< d_i$:
In this case we have $k_i\geq d_{i+1}-d_i>0$.
It suffices to prove that 
$$ a_ia_{i+1}^m a_{i-1}^{d_{i-1}}\ldots a_{i-1}^{k_{i-1}}a_i^m\ =\ a_{i+1}^m a_i a_{i-1}^{d_{i-1}}\ldots a_{i-1}^{k_{i-1}}a_i^m$$
for any $m>0$. Then the desired statement follows using equality~\eqref{power} from Lemma~\ref{lem:plactitech1}:
\begin{align*}
a_ia_{i+1}^{d_{i+1}}a_i^{d_i} a_{i-1}^{d_{i-1}}\ldots a_{i-1}^{k_{i-1}}a_i^{k_i}\ 
&=\ (a_{i+1}a_i)^{d_i}\left(a_ia_{i+1}^{d_{i+1}-d_i}a_{i-1}^{d_{i-1}}\ldots a_{i-1}^{k_{i-1}}a_i^{d_{i+1}-d_i}\right)a_i^{k_i+d_i-d_{i+1}}\\
&=\ (a_{i+1}a_i)^{d_i}\left(a_{i+1}^{d_{i+1}-d_i}a_i a_{i-1}^{d_{i-1}}\ldots a_{i-1}^{k_{i-1}}a_i^{d_{i+1}-d_i}\right)a_i^{k_i+d_i-d_{i+1}}\\
&=\ a_{i+1}^{d_{i+1}}a_i a_i^{d_i} a_{i-1}^{d_{i-1}}\ldots a_{i-1}^{k_{i-1}}a_i^{k_i}.
\end{align*}
Now for $d_{i-1},d_{i-2},\ldots,d_j>0$ and $d_{j-1}=0$ (possibly $j=i$, or $j=1$), we apply equations~\eqref{comm},~\eqref{part},~\eqref{power} and~\eqref{wurmcomm} to pass the factor $a_i$ (distinguished by bold print) through the whole expression, thereby proving the desired equality.
{\footnotesize
\begin{align*}
\ &\pmb{a_i}a_{i+1}^m a_{i-1}^{d_{i-1}}a_{i-2}^{d_{i-2}}\ldots a_j^{d_j}a_{j-2}^{d_{j-2}}\ldots a_{j-2}^{k_{j-2}}a_{j-1}^{k_{j-1}}a_j^{k_j}\ldots a_{i-1}^{k_{i-1}}a_i^m\\
\stackrel{\eqref{comm}}{=}\ &\pmb{a_i} a_{i-1}^{d_{i-1}}a_{i-2}^{d_{i-2}}\ldots a_j^{d_j}a_{j-2}^{d_{j-2}}\ldots a_{j-2}^{k_{j-2}}a_{j-1}^{k_{j-1}}a_j^{k_j}\ldots a_{i-1}^{k_{i-1}}a_{i+1}^m a_i^m\\
\stackrel{\text{\eqref{power}}}{=}\ &\pmb{a_i} a_{i-1}^{d_{i-1}}a_{i-2}^{d_{i-2}}\ldots a_j^{d_j}a_{j-2}^{d_{j-2}}\ldots a_{j-2}^{k_{j-2}}a_{j-1}^{k_{j-1}}a_j^{k_j}\ldots a_{i-1}^{k_{i-1}}(a_{i+1}a_i)^m\\
=\ &(\pmb{a_i} a_{i-1})a_{i-1}^{d_{i-1}-1}a_{i-2}^{d_{i-2}}\ldots a_j^{d_j}a_{j-2}^{d_{j-2}}\ldots a_{j-2}^{k_{j-2}}a_{j-1}^{k_{j-1}}a_j^{k_j}\ldots a_{i-1}^{k_{i-1}}(a_{i+1}a_i)^m\\
\stackrel{\text{\eqref{wurmcomm}}}{=}\ &a_{i-1}^{d_{i-1}-1}a_{i-2}^{d_{i-2}-1}\ldots (\pmb{a_i} a_{i-1} a_{i-2}\ldots a_j)a_j^{d_j-1}a_{j-2}^{d_{j-2}}\ldots a_{j-2}^{k_{j-2}}a_{j-1}^{k_{j-1}}a_j^{k_j}\ldots a_{i-1}^{k_{i-1}}(a_{i+1}a_i)^m\\
\stackrel{\text{\eqref{wurmcomm}}}{=}\ &a_{i-1}^{d_{i-1}-1}a_{i-2}^{d_{i-2}-1}\ldots a_j^{d_j-1}(\pmb{a_i} a_{i-1} a_{i-2}\ldots a_j)a_{j-2}^{d_{j-2}}\ldots a_{j-2}^{k_{j-2}}a_{j-1}^{k_{j-1}}a_j^{k_j}\ldots a_{i-1}^{k_{i-1}}(a_{i+1}a_i)^m\\
\stackrel{\eqref{comm}}{=}\ &a_{i-1}^{d_{i-1}-1}a_{i-2}^{d_{i-2}-1}\ldots a_j^{d_j-1}a_{j-2}^{d_{j-2}}\ldots a_{j-2}^{k_{j-2}}(\pmb{a_i} a_{i-1} a_{i-2}\ldots a_j)a_{j-1}^{k_{j-1}}a_j^{k_j}\ldots a_{i-1}^{k_{i-1}}(a_{i+1}a_i)^m\\
\stackrel{\text{\eqref{wurmcomm}}}{=}\ &a_{i-1}^{d_{i-1}-1}a_{i-2}^{d_{i-2}-1}\ldots a_j^{d_j-1}a_{j-2}^{d_{j-2}}\ldots a_{j-2}^{k_{j-2}}a_{j-1}^{k_{j-1}-1}(\pmb{a_i} a_{i-1} a_{i-2}\ldots a_j a_{j-1})a_j^{k_j}\ldots a_{i-1}^{k_{i-1}}(a_{i+1}a_i)^m\\
\stackrel{\text{\eqref{wurmcomm}}}{=}\ &a_{i-1}^{d_{i-1}-1}a_{i-2}^{d_{i-2}-1}\ldots a_j^{d_j-1}a_{j-2}^{d_{j-2}}\ldots a_{j-2}^{k_{j-2}}a_{j-1}^{k_{j-1}-1}a_j^{k_j}\ldots a_{i-1}^{k_{i-1}}(\pmb{a_i} a_{i-1} a_{i-2}\ldots a_j a_{j-1})(a_{i+1}a_i)^m\\
\stackrel{\eqref{comm}}{=}\ &a_{i-1}^{d_{i-1}-1}a_{i-2}^{d_{i-2}-1}\ldots a_j^{d_j-1}a_{j-2}^{d_{j-2}}\ldots a_{j-2}^{k_{j-2}}a_{j-1}^{k_{j-1}-1}a_j^{k_j}\ldots a_{i-1}^{k_{i-1}}(\pmb{a_i} a_{i-1}) (a_{i+1}a_i)^m(a_{i-2}\ldots a_j a_{j-1})\\
\stackrel{\eqref{part}}{=}\ &a_{i-1}^{d_{i-1}-1}a_{i-2}^{d_{i-2}-1}\ldots a_j^{d_j-1}a_{j-2}^{d_{j-2}}\ldots a_{j-2}^{k_{j-2}}a_{j-1}^{k_{j-1}-1}a_j^{k_j}\ldots a_{i-1}^{k_{i-1}}(a_{i+1}a_i)^m(\pmb{a_i} a_{i-1}) (a_{i-2}\ldots a_j a_{j-1})\\
\stackrel{\text{\eqref{power}}}{=}\ &a_{i-1}^{d_{i-1}-1}a_{i-2}^{d_{i-2}-1}\ldots a_j^{d_j-1}a_{j-2}^{d_{j-2}}\ldots a_{j-2}^{k_{j-2}}a_{j-1}^{k_{j-1}-1}a_j^{k_j}\ldots a_{i-1}^{k_{i-1}}a_{i+1}^m a_i^m(\pmb{a_i} a_{i-1}) (a_{i-2}\ldots a_j a_{j-1})\\
\stackrel{\eqref{comm}}{=}\ &a_{i+1}^m a_{i-1}^{d_{i-1}-1}a_{i-2}^{d_{i-2}-1}\ldots a_j^{d_j-1}a_{j-2}^{d_{j-2}}\ldots a_{j-2}^{k_{j-2}}a_{j-1}^{k_{j-1}-1}a_j^{k_j}\ldots a_{i-1}^{k_{i-1}} a_i^m(\pmb{a_i} a_{i-1} a_{i-2}\ldots a_j a_{j-1})\\
\stackrel{\text{\eqref{wurmcomm}}}{=}\ &a_{i+1}^m a_{i-1}^{d_{i-1}-1}a_{i-2}^{d_{i-2}-1}\ldots a_j^{d_j-1}a_{j-2}^{d_{j-2}}\ldots a_{j-2}^{k_{j-2}}(\pmb{a_i} a_{i-1} a_{i-2}\ldots a_j a_{j-1})a_{j-1}^{k_{j-1}-1}a_j^{k_j}\ldots a_{i-1}^{k_{i-1}} a_i^m\\
=\ &a_{i+1}^m a_{i-1}^{d_{i-1}-1}a_{i-2}^{d_{i-2}-1}\ldots a_j^{d_j-1}a_{j-2}^{d_{j-2}}\ldots a_{j-2}^{k_{j-2}}(\pmb{a_i} a_{i-1} a_{i-2}\ldots a_j)a_{j-1}^{k_{j-1}}a_j^{k_j}\ldots a_{i-1}^{k_{i-1}} a_i^m\\
\stackrel{\eqref{comm}}{=}\ &a_{i+1}^m a_{i-1}^{d_{i-1}-1}a_{i-2}^{d_{i-2}-1}\ldots a_j^{d_j-1}(\pmb{a_i} a_{i-1} a_{i-2}\ldots a_j)a_{j-2}^{d_{j-2}}\ldots a_{j-2}^{k_{j-2}}a_{j-1}^{k_{j-1}}a_j^{k_j}\ldots a_{i-1}^{k_{i-1}} a_i^m\\
\stackrel{\text{\eqref{wurmcomm}}}{=}\ &a_{i+1}^m a_{i-1}^{d_{i-1}-1}a_{i-2}^{d_{i-2}-1}\ldots (\pmb{a_i} a_{i-1} a_{i-2}\ldots a_j)a_j^{d_j-1}a_{j-2}^{d_{j-2}}\ldots a_{j-2}^{k_{j-2}}a_{j-1}^{k_{j-1}}a_j^{k_j}\ldots a_{i-1}^{k_{i-1}} a_i^m\\
\stackrel{\text{\eqref{wurmcomm}}}{=}\ &a_{i+1}^m a_{i-1}^{d_{i-1}-1}(\pmb{a_i} a_{i-1} a_{i-2})a_{i-2}^{d_{i-2}-1}\ldots a_j^{d_j}a_{j-2}^{d_{j-2}}\ldots a_{j-2}^{k_{j-2}}a_{j-1}^{k_{j-1}}a_j^{k_j}\ldots a_{i-1}^{k_{i-1}} a_i^m\\
\stackrel{\text{\eqref{wurmcomm}}}{=}\ &a_{i+1}^m (\pmb{a_i} a_{i-1}) a_{i-1}^{d_{i-1}-1}a_{i-2}^{d_{i-2}}\ldots a_j^{d_j}a_{j-2}^{d_{j-2}}\ldots a_{j-2}^{k_{j-2}}a_{j-1}^{k_{j-1}}a_j^{k_j}\ldots a_{i-1}^{k_{i-1}} a_i^m\\
=\ &a_{i+1}^m \pmb{a_i} a_{i-1}^{d_{i-1}}a_{i-2}^{d_{i-2}}\ldots a_j^{d_j}a_{j-2}^{d_{j-2}}\ldots a_{j-2}^{k_{j-2}}a_{j-1}^{k_{j-1}}a_j^{k_j}\ldots a_{i-1}^{k_{i-1}} a_i^m.
\end{align*}
}
\end{itemize}
\end{enumerate}
This concludes the proof of Proposition~\ref{prop:nofo}.
\end{bew}
Note that the relation special for the partic algebra~\eqref{part} was only used once in the proof of Proposition~\ref{prop:nofo}, namely in the long computation at the end. All other steps have been carried out using only the commutativity relation~\eqref{comm} and the plactic relations~\eqref{plac1} and~\eqref{plac2}. The following corollary recaps what we obtained for the multiplication in the partic algebra:
\begin{cor}\label{cor:particmulti}
Assume we are given a monomial $a_{N-1}^{d_{N-1}}\ldots a_2^{d_2} a_1^{k_1}a_2^{k_2}\ldots a_{N-1}^{k_{N-1}}$ of the form~\eqref{cond} in the partic algebra. Then left multiplication with $a_i$ gives
\begin{align}
&a_i\cdot\left(a_{N-1}^{d_{N-1}}\ldots a_2^{d_2} a_1^{k_1}a_2^{k_2}\ldots a_{N-1}^{k_{N-1}}\right)\label{leftmult}\\ 
=\ &\begin{cases}a_{N-1}^{d_{N-1}}\ldots a_i^{d_i}\ldots a_2^{d_2} a_1^{k_1}a_2^{k_2}\ldots a_i^{k_i+1}\ldots a_{N-1}^{k_{N-1}}&\text{ if }d_i=d_{i-1}+k_{i-1},\\
a_{N-1}^{d_{N-1}}\ldots a_i^{d_i+1}\ldots a_2^{d_2} a_1^{k_1}a_2^{k_2}\ldots a_i^{k_i}\ldots a_{N-1}^{k_{N-1}}&\text{ if }d_i<d_{i-1}+k_{i-1}.
\end{cases}\notag
\end{align}
Right multiplication with $a_i$ gives
\begin{align}
&\left(a_{N-1}^{d_{N-1}}\ldots a_2^{d_2} a_1^{k_1}a_2^{k_2}\ldots a_{N-1}^{k_{N-1}}\right)\cdot a_i\label{rightmult}\\
 =\ &\begin{cases}a_{N-1}^{d_{N-1}}\ldots a_i^{d_i}a_{i+1}^{d_{i+1}+1}\ldots a_2^{d_2} a_1^{k_1}a_2^{k_2}\ldots a_i^{k_i+1}a_{i+1}^{k_{i+1}-1}\ldots a_{N-1}^{k_{N-1}}&\text{ if }k_{i+1}\geq1\\
a_{N-1}^{d_{N-1}}\ldots a_i^{d_i}a_{i+1}^{d_{i+1}}\ldots a_2^{d_2} a_1^{k_1}a_2^{k_2}\ldots a_i^{k_i+1}a_{i+1}^0\ldots a_{N-1}^{k_{N-1}}&\text{ if }k_{i+1}=0,
\end{cases}\notag
\end{align}
with the result written again in the normal form~\eqref{cond}.
\end{cor}
\begin{bew}
For the left multiplication, equation~\eqref{leftmult} is contained in the proof of Proposition~\ref{prop:nofo}. For the right multiplication, equation~\eqref{rightmult} follows from the repeated application of the rule for left multiplication of $a_{N-1}^{d_{N-1}}\ldots a_2^{d_2} a_1^{k_1}a_2^{k_2}\ldots a_{N-1}^{k_{N-1}}$ to $a_i$.
\end{bew}
\begin{bsp}\label{ex:nottfree}
The partic algebra has zero divisors, e.g. in $\partic{3}$,
$$a_2\cdot\left(a_3^5a_2^8a_1^8a_2^3a_3^1\ -\ a_3^5a_2^7a_1^8a_2^4a_3^1\right)\ =\ a_3^5a_2^8a_1^8a_2^4a_3^1\ -\ a_3^5a_2^8a_1^8a_2^4a_3^1\ =\ 0$$
(it follows from Theorem~\ref{theo:particbasis} that $a_3^5a_2^8a_1^8a_2^3a_3^1\ -\ a_3^5a_2^7a_1^8a_2^4a_3^1\ \neq\ 0$).
\end{bsp}
\begin{bem}
Let us compare our normal form with the monomial bases of the plactic algebra from \cite{reinmono}: The plactic algebra $\plac{N}$ surjects onto the partic algebra $\partic{N}$, mapping generators to generators and hence monomials to monomials. Given a monomial of the normal form from Proposition~\ref{prop:nofo}, finding the (finitely many) preimages of basis monomials in the plactic algebra amounts to solving a system of linear equations over the nonnegative integers, i.e. finding lattice points inside a polyhedron.
\par
For example, consider the basis of the plactic algebra $\plac{5}$ from \cite[Theorem~2.10]{reinmono} given by monomials
$$\resizebox{.9\hsize}{!}{$(a_1)^{n_1}(a_2a_1)^{n_{21}}(a_2)^{n_2}(a_3a_2a_1)^{n_{321}}(a_3a_2)^{n_{32}}(a_3)^{n_{3}}(a_4a_3a_2a_1)^{n_{4321}}(a_4a_3a_2)^{n_{432}}(a_4a_3)^{n_{43}}(a_4)^{n_4}$}$$
where all $n_i\in\ZZ_{\geq0}$ and compare it with the basis of the partic algebra $\partic{5}$ from Proposition~\ref{prop:nofo}
$$\{ a_4^{d_4}a_3^{d_3}a_2^{d_2}a_1^{k_1}a_2^{k_2}a_3^{k_3}a_4^{k_4}\ |\ \text{all }k_i,d_i\in\ZZ_{\geq0},\ d_i\leq d_{i-1}+k_{i-1}\}.$$
While $a_1a_2a_3a_4\in\partic{5}$ has only one preimage, namely $(a_1)^1(a_2)^1(a_3)^1(a_4)^1\in\plac{5}$, we find two preimages of $a_4a_3a_2a_1a_2\in\partic{5}$, namely $(a_2)^1(a_4a_3a_2a_1)^1,(a_2a_1)^1(a_4a_3a_2)^1\in\plac{5}$. This corresponds to the number of possible applications of the additional partic relation~\eqref{part}. 
\end{bem}

\section{The action on bosonic particle configurations}\label{sec:bosonicac}
In this section we discuss an action of the plactic algebra $\plac{N}$ on the polynomial ring $\ground[x_1,\ldots,x_{N-1},x_0]$ in $N$ variables. It was defined in \cite[Proposition~5.8]{ks}. We recall the definition here:
Let $x_1^{k_1}\ldots x_{N-1}^{k_{N-1}}x_0^{k_0}$ be a monomial in $\ground[x_1,\ldots,x_{N-1},x_0]$. Set
\begin{align}\label{eq:placticaction}
 a_i\cdot x_1^{k_1}\ldots x_{N-1}^{k_{N-1}}x_0^{k_0}\ &=\ \begin{cases} x_1^{k_1}\ldots x_i^{k_i-1}x_{i+1}^{k_{i+1}+1}\ldots x_{N-1}^{k_{N-1}}x_0^{k_0}\quad &\text{if }k_i>0,\\ 0&\text{else,}\end{cases}\\
a_{N-1}\cdot x_1^{k_1}\ldots x_{N-1}^{k_{N-1}}x_0^{k_0}\ &=\ \begin{cases} x_1^{k_1}\ldots  x_{N-1}^{k_{N-1}-1}x_0^{k_0+1}\quad &\text{if }k_{N-1}>0,\\ 0&\text{else.}\end{cases}\label{eq:placticactionN}
\end{align}
This defines an action of the plactic algebra which factors over the partic algebra:
\begin{lemma}
Equations~\eqref{eq:placticaction} and~\eqref{eq:placticactionN} define an action of the plactic algebra $\plac{N}$ on the polynomial ring $\ground[x_1,\ldots,x_{N-1},x_0]$. This action factors over an action of the partic algebra $\partic{N}$. 
\end{lemma}
\begin{bew}
This can be verified by direct computation.
\end{bew}
In this section our goal is the proof of the following main theorem:
\begin{theo}\label{theo:particfaith}
The action of the partic algebra $\partic{N}$ on $\ground[x_1,\ldots,x_{N-1},x_0]$ defined by equations~\eqref{eq:placticaction} and~\eqref{eq:placticactionN} is faithful.
\end{theo}
\begin{bem}
In \cite[Proposition~5.8]{ks} it is stated incorrectly that the action of the plactic algebra $\plac{N}$ on $\ground[x_1,\ldots,x_{N-1},x_0]$ is faithful.
\end{bem}
\begin{defi}\label{defi:partinot}
We introduce the shorthand notation $\bosconfi{i}:=(k_1,\ldots,k_{N-1},k_0)\in\ZZ_{\geq0}^{N}$ for the monomial $v(\bosconfi{i}):=x_1^{k_1}\ldots x_{N-1}^{k_{N-1}}x_0^{k_0}$.
\end{defi}
One can think of the monomial $x_1^{k_1}\ldots x_{N-1}^{k_{N-1}}x_0^{k_0}$ or the tuple $(k_1,\ldots,k_{N-1},k_0)$ as a configuration of particles on a line with $N$ positions, with $k_i$ particles at the $i$-th position. The $0$-th position is regarded as the deposit for particles moved to the end of the line. Then $a_i$ moves a particle from position $i$ to position $i+1$. We call $\ground[x_1,\ldots,x_{N-1},x_0]$ with the above action the (classical bosonic) particle configuration module of $\plac{N}$ or $\partic{N}$, and we refer to the monomials inside $\ground[x_1,\ldots,x_{N-1},x_0]$ as (classical bosonic) particle configurations. 

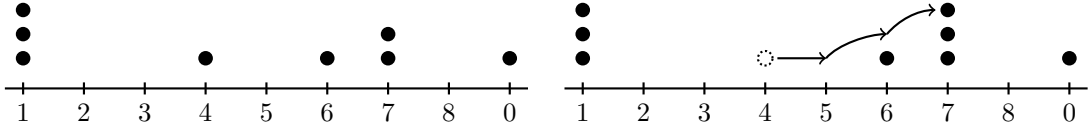
\begin{figure}[H]
\begin{center}
\begin{tikzpicture}[scale=0.8]

\begin{scope}[xshift=-4.6cm]
\path[use as bounding box] (-4.5,-1.5) rectangle (4.5,2.5);
\draw[thick] (-4.3,0) -- (4.3,0);

\foreach \x in {1,...,8}
\draw[thick] ({-5+\x},-0.1) -- ({-5+\x},0.1);

\draw[thick] (4,-0.1) -- (4,0.1);

\foreach \x in {1,...,8}
\draw ({-5+\x},-0.4) node {$\x$};

\draw (4,-0.4) node {$0$};

\foreach \x in {1,4,6,7,9}
\fill ({-5+\x},0.5) circle (1.2mm);
\foreach \x in {1,7}
\fill ({-5+\x},0.9) circle (1.2mm);
\foreach \x in {1}
\fill ({-5+\x},1.3) circle (1.2mm);

\end{scope}
\begin{scope}[xshift=4.6cm]
\path[use as bounding box] (-4.5,-1.5) rectangle (4.5,2.5);
\draw[thick] (-4.3,0) -- (4.3,0);

\foreach \x in {1,...,8}
\draw[thick] ({-5+\x},-0.1) -- ({-5+\x},0.1);

\draw[thick] (4,-0.1) -- (4,0.1);

\foreach \x in {1,...,8}
\draw ({-5+\x},-0.4) node {$\x$};

\draw (4,-0.4) node {$0$};

\foreach \x in {1,6,7,9}
\fill ({-5+\x},0.5) circle (1.2mm);
\foreach \x in {1,7}
\fill ({-5+\x},0.9) circle (1.2mm);
\foreach \x in {1,7}
\fill ({-5+\x},1.3) circle (1.2mm);

\foreach \x in {4}
\draw[thick, densely dotted] ({-5+\x},0.5) circle (1.2mm);

\draw[thick, ->] ({-0.8},{0.5}) -- ({0},{0.5});
\draw[thick, ->] ({0},{0.5}) .. controls +(50:0.3cm) and +(180:0.3cm) .. ({1},{0.9});
\draw[thick, ->] ({1},{0.9}) .. controls +(50:0.3cm) and +(180:0.3cm) .. (1.8,1.3);

\end{scope}
\end{tikzpicture}
\end{center}
\caption{Example for $N=9$: The particle configuration $(3,0,0,1,0,1,2,0,1)$ corresponding to the monomial 
$x_1^3 x_2^0 x_3^0 x_4^1 x_5^0 x_6^1 x_7^2 x_8^0 x_0^1$, and the element $a_6 a_5 a_4$ acting on it.}
\label{fig:bosonicclass}
\end{figure}
\par
Now we investigate the action of the partic algebra on the particle configuration module.
\begin{prop}\label{prop:labelmonomialpartic}
Fix a monomial
$a_{N-1}^{d_{N-1}}\ldots a_3^{d_3}a_2^{d_2} a_1^{k_1}a_2^{k_2}a_3^{k_3}\ldots a_{N-1}^{k_{N-1}}$ in the partic algebra satisfying condition~\eqref{cond}. There is a unique particle configuration with the number of particles minimal, i.e. a monomial in $\ground[x_1,\ldots,x_{N-1},x_0]$ of minimal degree, so that the monomial acts nontrivially on it. This minimal particle configuration is given by
$$\bosconfi{i}_\text{in}\ = \ (k_1,k_2,k_3\ldots,k_{N-1},0).$$
The image of $\bosconfi{i}_\text{in}$ under the action of $a_{N-1}^{d_{N-1}}\ldots a_2^{d_2} a_1^{k_1}a_2^{k_2}\ldots a_{N-1}^{k_{N-1}}$ is given by
$$\bosconfi{i}_\text{out}\ =\ (0,k_1-d_2,k_2+d_2-d_3\ldots, k_{N-2}+d_{N-2}-d_{N-1}, k_{N-1}+d_{N-1}).$$
\end{prop}
\begin{bew}
First we show that $a_1^{k_1}a_2^{k_2}a_3^{k_3}\ldots a_{N-1}^{k_{N-1}}$, hence $a_{N-1}^{d_{N-1}}\ldots a_3^{d_3}a_2^{d_2} a_1^{k_1}a_2^{k_2}a_3^{k_3}\ldots a_{N-1}^{k_{N-1}}$ annihilates any particle configuration $(r_1,r_2,r_3,\ldots,r_{N-1},r_0)$ with $r_i<k_i$ for some $i$.
We compute
{\footnotesize
\begin{align*}
&a_1^{k_1}a_2^{k_2}\ldots a_j^{k_j}a_{j+1}^{k_{j+1}} \ldots a_{N-1}^{k_{N-1}}(x_1^{r_1}x_2^{r_2}\ldots x_j^{r_j}x_{j+1}^{r_{j+1}}\ldots x_{N-1}^{r_{N-1}}x_0^{r_0})\\
=\ &\begin{cases}a_1^{k_1}a_2^{k_2}\ldots a_j^{k_j}(x_1^{r_1}x_2^{r_2}\ldots x_j^{r_j}x_{j+1}^{r_{j+1}-k_{j+1}}\ldots x_{N-1}^{r_{N-1}-k_{N-1}+k_{N-2}}x_0^{r_0+k_{N-1}}), &\text{}r_{i}\geq k_{i}\text{ for }j < i\leq {N-1}\\ 0&\text{else}\end{cases}\\
=\ &\begin{cases}x_1^{r_1-k_1}x_2^{r_2-k_2+k_1}\ldots x_i^{r_i-k_i+k_{i-1}}x_{i+1}^{r_{i+1}-k_{i+1}+k_i}\ldots x_{N-1}^{r_{N-1}-k_{N-1}+k_{N-2}}x_0^{r_0+k_{N-1}},&\text{}r_i\geq k_i\text{ for all }i\\ 0&\text{else.}\end{cases}
\end{align*}
}
Together with condition~\eqref{cond} it follows that the action of a monomial of the form 
$$a_{N-1}^{d_{N-1}}\ldots a_3^{d_3}a_2^{d_2} a_1^{k_1}a_2^{k_2}a_3^{k_3}\ldots a_{N-1}^{k_{N-1}}$$
on a particle configuration $(r_1,r_2,r_3,\ldots,r_{N-1},r_0)$ is nontrivial iff $r_i\geq k_i$ for all $i$ (recall that $r_0\geq0=k_0$ is automatically satisfied). This proves that $\bosconfi{i}_\text{in}$ is indeed the minimal particle configuration on which the monomial acts nontrivially. Now compute the image of $\bosconfi{i}_\text{in}$ under the action of the monomial: Plug in $r_i=k_i$ for all $i$ to see that
\begin{align*}
a_{N-1}^{d_{N-1}}\ldots a_3^{d_3}a_2^{d_2} a_1^{k_1}a_2^{k_2}a_3^{k_3}\ldots a_{N-1}^{k_{N-1}}(\bosconfi{i}_\text{in})\ &=\ a_{N-1}^{d_{N-1}}\ldots a_3^{d_3}a_2^{d_2}(x_1^{0}x_2^{k_1}\ldots x_{N-1}^{k_{N-2}}x_0^{k_{N-1}})\\
&=\ \bosconfi{i}_\text{out}.
\end{align*}
This proves Proposition~\ref{prop:labelmonomialpartic}.
\end{bew}
\begin{bew}[of Theorem~\ref{theo:particbasis}]
By Proposition~\ref{prop:nofo} any monomial in the partic algebra is equivalent to one of the form~\eqref{cond}. We have shown in Proposition~\ref{prop:labelmonomialpartic} that the action on the particle configuration module distinguishes any two monomials of the form~\eqref{cond}, hence~\eqref{cond} describes a normal form for the monomials in the partic algebra $\partic{N}$, hence a basis of $\partic{N}$.
\end{bew}
Now Theorem~\ref{theo:particfaith} follows as a corollary from Proposition~\ref{prop:labelmonomialpartic}:
\begin{bew}[of Theorem~\ref{theo:particfaith}]
We have seen in Proposition~\ref{prop:labelmonomialpartic} that the normal form monomials, hence the basis elements in $\partic{N}$ act linearly independent on the particle configurations. In other words, the action of $\partic{N}$ is faithful.
\end{bew}
\begin{bem}
The faithfulness of the action of the algebra $\partic{N}$ on the particle configuration module motivates us to give $\partic{N}$ the name ``partic'' algebra.
\end{bem}
\par
By Proposition~\ref{prop:nofo} and Proposition~\ref{prop:labelmonomialpartic}, we can identify each monomial in the partic algebra uniquely by the minimal particle configuration $\bosconfi{j}\in\ZZ_{\geq0}^{N}$ on which it acts nontrivially and the output particle configuration $\bosconfi{i}\in\ZZ_{\geq0}^{N}$ that one gets back from the action of the monomial on $\bosconfi{j}$. Hence the following is welldefined:
\begin{defi}\label{defi:particaij}
Given a monomial in normal form with $d_i\leq d_{i-1}+k_{i-1}$ for all $3\leq i\leq {N-1}$ and $d_2\leq k_1$, see Proposition~\ref{prop:labelmonomialpartic}, we write
$$a_{\bosconfi{i}\bosconfi{j}}=a_{N-1}^{d_{N-1}}\ldots a_2^{d_2} a_1^{k_1}a_2^{k_2}\ldots a_{N-1}^{k_{N-1}}$$
for bosonic particle configurations $\bosconfi{i}=(0,k_1-d_2,k_2+d_2-d_3\ldots, k_{N-2}+d_{N-2}-d_{{N-1}}, k_{N-1}+d_{N-1})$ and $\bosconfi{j}=(k_1,k_2,k_3\ldots,k_{N-1},0)$.
The number of particles $|\bosconfi{i}|=|\bosconfi{j}|=\sum_i k_i$ in $\bosconfi{i}$ and $\bosconfi{j}$ is the same.
\end{defi}
This labelling is made so that $a_{\bosconfi{i}\bosconfi{j}}\cdot v(\bosconfi{j})=v(\bosconfi{i})$ in the notation of Definition \ref{defi:partinot}.
\begin{defi}
For $\bosconfi{i}=(r_1,\ldots,r_{N-1},r_0)\in\ZZ_{\geq0}^{N}$, we set 
$$\bosconfi{i}\cup\{i\} = (r_1,\ldots,r_i+1,\ldots,r_{N-1},r_0),\qquad \bosconfi{i}\setminus\{i\} = (r_1,\ldots,r_i-1,\ldots,r_{N-1},r_0),$$
where the latter is only defined for $r_i>0$.
\end{defi}
With this notation we can rewrite Corollary \ref{cor:particmulti} to obtain the following multiplication rule.
\begin{cor}\label{cor:particmulticonfi}
Let $a_{\bosconfi{i}\bosconfi{j}}$ be a monomial in normal form as in Definition \ref{defi:particaij}. Then left and right multiplication by some generator $a_i\in\partic{N}$ are given by
\begin{align*}
a_i a_{\bosconfi{i}\bosconfi{j}}\ &=\ \begin{cases} a_{\bosconfi{i}\pp\bosconfi{j}\pp}\quad &\text{if }i\in\bosconfi{i}\\ a_{{\bosconfi{i}}\p{\bosconfi{j}}\p}\quad &\text{if }i\notin\bosconfi{i},\end{cases}\\
a_{\bosconfi{i}\bosconfi{j}}a_i \ &=\ \begin{cases} a_{{\bosconfi{i}\ppp}{\bosconfi{j}\ppp}}\quad &\text{if }i+1\in\bosconfi{j}\\ a_{\bosconfi{i}\p\bosconfi{j}\p}\quad &\text{if }i+1\notin\bosconfi{j}.\end{cases}
\end{align*}
Here we denote
\begin{align*}
\bosconfi{i}\p\ &=\ \bosconfi{i}\cup\{i+1\}\quad &\bosconfi{i}\pp\ &=\ (\bosconfi{i}\setminus\{i\})\cup\{i+1\}\quad &{\bosconfi{i}\ppp}\ &=\ \bosconfi{i}\\
\bosconfi{j}\p\ &=\ \bosconfi{j}\cup\{i\}\quad &\bosconfi{j}\pp\ &=\ \bosconfi{j}\quad &{\bosconfi{j}\ppp}\ &=\ (\bosconfi{j}\setminus\{i+1\})\cup\{i\}.
\end{align*}
\end{cor}
\begin{bsp}
Let $N=6$, and consider the monomial $a_{\bosconfi{i}\bosconfi{j}}=a_5^1a_2^2a_3^1a_4^2$ with minimal input configurataion $\bosconfi{j}=(0,2,1,2,0,0)$ and output configuration $\bosconfi{i}=(0,0,2,1,1,1)$. Now consider the left and right multiplication with $a_i$ for $i=3$: 
\begin{align*}
a_3 \cdot a_{(0,0,2,1,1,1)(0,2,1,2,0,0)}\ &=\ a_{(0,0,1,2,1,1)(0,2,1,2,0,0)},\\
\text{with}\quad\bosconfi{i}\pp\ &=\ (0,0,1,2,1,1),\\
\bosconfi{j}\pp\ &=\ (0,2,1,2,0,0),\\
a_{(0,0,2,1,1,1)(0,2,1,2,0,0)}a_3 \ &=\ a_{\bosconfi(0,0,2,1,1,1)(0,2,2,1,0,0)},\\
\text{with}\quad\bosconfi{i}\ppp\ &=\ (0,0,2,1,1,1),\\
\bosconfi{j}\ppp\ &=\ (0,2,2,1,0,0).
\end{align*}
In contrast, left and right multiplication with $a_i$ for $i=1$ gives
\begin{align*}
a_1 \cdot a_{(0,0,2,1,1,1)(0,2,1,2,0,0)}\ &=\ a_{(0,1,1,2,1,1)(1,2,1,2,0,0)},\\
\text{with}\quad\bosconfi{i}\p\ &=\ (0,1,1,2,1,1),\\
\bosconfi{j}\p\ &=\ (1,2,1,2,0,0),\\ 
a_{(0,0,2,1,1,1)(0,2,1,2,0,0)}a_1 \ &=\ a_{(0,0,2,1,1,1)(1,1,1,2,0,0)},\\
\text{with}\quad\bosconfi{i}\ppp\ &=\ (0,0,2,1,1,1),\\
\bosconfi{j}\ppp\ &=\ (1,1,1,2,0,0).
\end{align*}
We observe that the product $a_1 \cdot a_{(0,0,2,1,1,1)(0,2,1,2,0,0)}$ requires an additional particle at position $1$, so that the cardinality of the minimal particle configuration of the product $a_1 \cdot a_{(0,0,2,1,1,1)(0,2,1,2,0,0)}$ is by one higher than that of $a_{(0,0,2,1,1,1)(0,2,1,2,0,0)}$.
\end{bsp}
\section{The center of the partic algebra}\label{sec:bosce}
Now that we have a basis of the partic algebra with a convenient labelling at our disposal, the goal of this section is to describe the center of the partic algebra $\partic{N}$.
\begin{theo}\label{theo:centpart}
The center of the partic algebra $\partic{N}$ is given by the $\ground$-span of the elements
$$\{ a_{N-1}^r a_{N-2}^r\ldots a_2^r a_1^r\ |\ r\geq0\}.$$
\end{theo}
The monomial $a_{N-1}^r a_{N-2}^r\ldots a_2^r a_1^r=(a_{N-1} a_{N-2}\ldots a_2 a_1)^r = a_{(0,\ldots,0,r)(r,0,\ldots,0)}$ acts on the bosonic particle configurations by moving $r$ particles from the first position $1$ to the last position $0$ if there are at least $r$ particles at position $1$, and it acts by zero if there are less than $r$ particles at position $1$. This action can be visualized as follows:
\begin{figure}[H]
\begin{center}
\begin{tikzpicture}[scale=0.8]

\begin{scope}[xshift=0cm]
\path[use as bounding box] (-4.5,-1.5) rectangle (4.5,2.5);
\draw[thick] (-4.3,0) -- (4.3,0);

\foreach \x in {1,...,8}
\draw[thick] ({-5+\x},-0.1) -- ({-5+\x},0.1);

\draw[thick] (4,-0.1) -- (4,0.1);

\foreach \x in {1,...,8}
\draw ({-5+\x},-0.4) node {$\x$};

\draw (4,-0.4) node {$0$};

\foreach \x in {1}
\draw[thick, densely dotted] ({-5+\x},0.5) circle (1.2mm);
\foreach \x in {1}
\draw[thick, densely dotted] ({-5+\x},0.9) circle (1.2mm);
\foreach \x in {1}
\draw[thick, densely dotted] ({-5+\x},1.3) circle (1.2mm);
\foreach \x in {1}
\draw[thick, densely dotted] ({-5+\x},1.7) circle (1.2mm);
\foreach \x in {1}
\draw[thick, densely dotted] ({-5+\x},2.1) circle (1.2mm);

\foreach \x in {9}
\fill ({-5+\x},0.5) circle (1.2mm);
\foreach \x in {9}
\fill ({-5+\x},0.9) circle (1.2mm);
\foreach \x in {9}
\fill ({-5+\x},1.3) circle (1.2mm);
\foreach \x in {9}
\fill ({-5+\x},1.7) circle (1.2mm);
\foreach \x in {9}
\fill ({-5+\x},2.1) circle (1.2mm);

\foreach \x in {2,...,7}
\foreach \y in {0,...,4}
\draw[thick, ->] ({-5+\x},{0.5+\y*0.4}) -- ({-4+\x},{0.5+\y*0.4});

\foreach \y in {0,...,4}
\draw[thick, ->] ({-3.8},{0.5+\y*0.4}) -- ({-3},{0.5+\y*0.4});

\foreach \y in {0,...,4}
\draw[thick, ->] ({3},{0.5+\y*0.4}) -- ({3.8},{0.5+\y*0.4});
\end{scope}
\end{tikzpicture}
\end{center}
\caption{Example for $N=9$: The action of the central element $(a_8a_7a_6a_5a_4a_3a_2a_1)^5 $ on the particle configuration $(5,0,0,0,0,0,0,0,0)$.}
\end{figure}
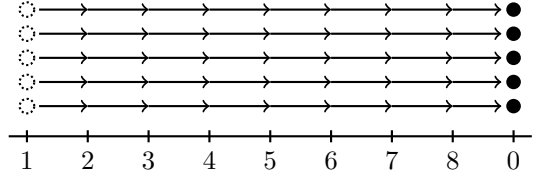
\begin{bew}
Let $z:=\sum\limits_{\bosconfi{i},\bosconfi{j}} c_{\bosconfi{i}\bosconfi{j}} a_{\bosconfi{i}\bosconfi{j}}$ be an element in the center, where we label the monomial $a_{\bosconfi{i}\bosconfi{j}}$ by minimal input and output particle configurations as in Definition \ref{defi:particaij}, with coefficients $c_{\bosconfi{i}\bosconfi{j}}\in\ground$.
Notice that $a_{(0,\ldots,0,r)(r,0,\ldots,0)}$ commutes with all $a_i$ by equation~\eqref{wurmcomm} from Lemma~\ref{lem:plactitech2}. We show that $c_{\bosconfi{i}\bosconfi{j}}=0$ for all $\bosconfi{j}$ that contain some $i\neq1$, and for all $\bosconfi{i}$ that contain some $i\neq 0$.
\par
Let $i\geq 2$. First we prove that $c_{\bosconfi{i}\bosconfi{j}}=0$ for all $\bosconfi{j}$ that contain a particle at position $i$. Since $z=\sum\limits_{\bosconfi{i},\bosconfi{j}} c_{\bosconfi{i}\bosconfi{j}} a_{\bosconfi{i}\bosconfi{j}}$ is central, it commutes in particular with $a_{i-1}a_{i-2}\ldots a_2 a_1$. Using Corollary \ref{cor:particmulticonfi} we calculate
\begin{align*}
(a_{i-1}a_{i-2}\ldots a_2 a_1) a_{\bosconfi{i}\bosconfi{j}}\ &=\ a_{(\bosconfi{i}\cup\{i\})(\bosconfi{j}\cup\{1\})},\\
a_{\bosconfi{i}\bosconfi{j}}(a_{i-1}a_{i-2}\ldots a_2 a_1) \ &=\ \begin{cases}
a_{(\bosconfi{i}\cup\{i\})(\bosconfi{j}\cup\{1\})}\ &\text {if }i\notin\bosconfi{j},\\
a_{\bosconfi{i}((\bosconfi{j}\setminus\{i\})\cup\{1\})}\ &\text {if }i\in\bosconfi{j}.
\end{cases}
\end{align*}
Therefore $(a_{i-1}a_{i-2}\ldots a_2 a_1) a_{\bosconfi{i}\bosconfi{j}}= a_{\bosconfi{i}\bosconfi{j}}(a_{i-1}a_{i-2}\ldots a_2 a_1)$ for $i\notin\bosconfi{j}$. This we use to deduce that we have $(a_{i-1}a_{i-2}\ldots a_2 a_1)z=z(a_{i-1}a_{i-2}\ldots a_2 a_1)$ if and only if
$$
(a_{i-1}a_{i-2}\ldots a_2 a_1)\left(\sum\limits_{\substack{\bosconfi{i},\ \bosconfi{j}\\ i\notin\bosconfi{j}}} c_{\bosconfi{i}\bosconfi{j}} a_{\bosconfi{i}\bosconfi{j}}\ +\ \sum\limits_{\substack{\bosconfi{i},\ \bosconfi{j}\\ i\in\bosconfi{j}}} c_{\bosconfi{i}\bosconfi{j}} a_{\bosconfi{i}\bosconfi{j}}\right)\ =\ \left(\sum\limits_{\substack{\bosconfi{i},\ \bosconfi{j}\\ i\notin\bosconfi{j}}} c_{\bosconfi{i}\bosconfi{j}} a_{\bosconfi{i}\bosconfi{j}}\ +\ \sum\limits_{\substack{\bosconfi{i},\ \bosconfi{j}\\ i\in\bosconfi{j}}} c_{\bosconfi{i}\bosconfi{j}} a_{\bosconfi{i}\bosconfi{j}}\right)(a_{i-1}a_{i-2}\ldots a_2 a_1),
$$
which holds if and only if
\begin{align*}(a_{i-1}a_{i-2}\ldots a_2 a_1)\left(\sum\limits_{\substack{\bosconfi{i},\ \bosconfi{j}\\ i\in\bosconfi{j}}} c_{\bosconfi{i}\bosconfi{j}} a_{\bosconfi{i}\bosconfi{j}}\right)\ =\ \left(\sum\limits_{\substack{\bosconfi{i},\ \bosconfi{j}\\ i\in\bosconfi{j}}} c_{\bosconfi{i}\bosconfi{j}} a_{\bosconfi{i}\bosconfi{j}}\right)(a_{i-1}a_{i-2}\ldots a_2 a_1).
\end{align*}
The latter is precisely the equality
\begin{align}\label{leftsumbos}
&\sum\limits_{\substack{\bosconfi{i},\ \bosconfi{j}\\ i\in\bosconfi{j}}} c_{\bosconfi{i}\bosconfi{j}}a_{(\bosconfi{i}\cup\{i\})(\bosconfi{j}\cup\{1\})}
\ =\ 
\sum\limits_{\substack{\bosconfi{i},\ \bosconfi{j}\\ i\in\bosconfi{j}}} c_{\bosconfi{i}\bosconfi{j}}a_{\bosconfi{i}((\bosconfi{j}\setminus\{i\})\cup\{1\})}.
\end{align}
Observe on the other hand that for fixed $i$ the set of monomials
$$\{ a_{\bosconfi{i}((\bosconfi{j}\setminus\{i\})\cup\{1\})}\ |\ \bosconfi{i},\ \bosconfi{j}\text{ such that }i\in\bosconfi{j}\}$$
is linearly independent since the sets $((\bosconfi{j}\setminus\{i\})\cup\{1\})$ are all distinct for distinct $\bosconfi{j}$.
\par
Next, we show by induction on the number $k_i$ of particles at position $i$ in $\bosconfi{j}$ that all coefficients  $c_{\bosconfi{i}\bosconfi{j}}$ are zero for $k_i\geq1$:\\
For $k_i=1$, the set $(\bosconfi{j}\setminus\{i\})\cup\{1\}$ does not contain any particle at position $i$ any more. Hence the monomial $a_{\bosconfi{i}((\bosconfi{j}\setminus\{i\})\cup\{1\})}$ cannot appear in the left sum in equation~\eqref{leftsumbos}, and so its coefficient $c_{\bosconfi{i}\bosconfi{j}}$ must have been zero.
For the induction step, assume that the coefficient $c_{\bosconfi{i}\bosconfi{j}}$ is zero for all $a_{\bosconfi{i}\bosconfi{j}}$ with at most $k_i$ particles at position $i$ in the minimal input particle configuration $\bosconfi{j}$. Consider some $a_{\bosconfi{i}\bosconfi{j}}$ with $k_i+1$ particles at position $i$ in $\bosconfi{j}$.
So the set $(\bosconfi{j}\setminus\{i\})\cup\{1\}$ contains $k_i$ particles at position $i$ in $\bosconfi{j}$, and so the monomial $a_{\bosconfi{i}((\bosconfi{j}\setminus\{i\})\cup\{1\})}$ cannot appear in the sum~\eqref{leftsumbos}. Therefore we see that the coefficient $c_{\bosconfi{i}\bosconfi{j}}$ must have been zero.
\par
We have shown that any central element in $\partic{N}$ is of the form 
$$z\ =\ \sum\limits_{\bosconfi{i},\ \bosconfi{j}} c_{\bosconfi{i}\bosconfi{j}} a_{\bosconfi{i}\bosconfi{j}},$$
where the particle configurations $\bosconfi{j}$ are of the form $(r,0,\ldots,0)$, $r\in\ZZ_{\geq0}$.
We use the convention that $i+1=0$ for $i={N-1}$ which matches our definition of the action of the partic algebra $\partic{N}$ on the bosonic particle configuration module. Notice that $0$ is never contained in the minimal input particle configuration, so that for  $1\leq i\leq {N-1}$ we have that $i+1\notin\bosconfi{j}$ for all $c_{\bosconfi{i}\bosconfi{j}}\neq0$.
\par
Now we use a similar induction argument to show that $c_{\bosconfi{i}\bosconfi{j}}=0$ for all $\bosconfi{i}$ that contain a particle at position $i\neq0$.
So let $1\leq i\leq {N-1}$. 
Using Corollary \ref{cor:particmulticonfi} we calculate that 
\begin{align*}
a_i a_{\bosconfi{i}\bosconfi{j}}\ &=\ \begin{cases}
a_{(\bosconfi{i}\cup\{i+1\})(\bosconfi{j}\cup\{i\})}\ &\text {if }i\notin\bosconfi{i},\\
a_{((\bosconfi{i}\setminus\{i\})\cup\{i+1\})\bosconfi{j}}\ &\text {if }i\in\bosconfi{i},
\end{cases}\\
a_{\bosconfi{i}\bosconfi{j}}a_i \ &=\ \begin{cases}
a_{(\bosconfi{i}\cup\{i+1\})(\bosconfi{j}\cup\{i\})}\ &\text {if }i+1\notin\bosconfi{j},\\
a_{\bosconfi{i}((\bosconfi{j}\setminus\{i+1\})\cup\{i\})}\ &\text {if }i+1\in\bosconfi{j}.
\end{cases}
\end{align*}
Since we have shown already that $i+1\notin\bosconfi{j}$, we know that $a_i z=z a_i$ is nothing but the equality
$$a_i\left(\sum\limits_{\substack{\bosconfi{i},\ \bosconfi{j}\\ i+1\notin\bosconfi{j}\\ i\notin\bosconfi{i}}} c_{\bosconfi{i}\bosconfi{j}} a_{\bosconfi{i}\bosconfi{j}}\ +\ \sum\limits_{\substack{\bosconfi{i},\ \bosconfi{j}\\ i+1\notin\bosconfi{j}\\ i\in\bosconfi{i}}} c_{\bosconfi{i}\bosconfi{j}} a_{\bosconfi{i}\bosconfi{j}}\right)
\ =\ \left(\sum\limits_{\substack{\bosconfi{i},\ \bosconfi{j}\\ i+1\notin\bosconfi{j}\\ i\notin\bosconfi{i}}} c_{\bosconfi{i}\bosconfi{j}} a_{\bosconfi{i}\bosconfi{j}}\ +\ \sum\limits_{\substack{\bosconfi{i},\ \bosconfi{j}\\ i+1\notin\bosconfi{j}\\ i\in\bosconfi{i}}} c_{\bosconfi{i}\bosconfi{j}} a_{\bosconfi{i}\bosconfi{j}}\right)a_i.$$
This in turn is equivalent to the equality
\begin{align*}\notag
&a_i\left(\sum\limits_{\substack{\bosconfi{i},\ \bosconfi{j}\\ i+1\notin\bosconfi{j}\\ i\in\bosconfi{i}}} c_{\bosconfi{i}\bosconfi{j}} a_{\bosconfi{i}\bosconfi{j}}\right)\ =\ 
\left(\sum\limits_{\substack{\bosconfi{i},\ \bosconfi{j}\\ i+1\notin\bosconfi{j}\\ i\in\bosconfi{i}}} c_{\bosconfi{i}\bosconfi{j}} a_{\bosconfi{i}\bosconfi{j}}\right)a_i,
\end{align*}
which can be rewritten as
\begin{align}
\sum\limits_{\substack{\bosconfi{i},\ \bosconfi{j}\\ i+1\notin\bosconfi{j}\\ i\in\bosconfi{i}}} c_{\bosconfi{i}\bosconfi{j}} a_{((\bosconfi{i}\setminus\{i\})\cup\{i+1\})\bosconfi{j}}
\ =\ \sum\limits_{\substack{\bosconfi{i},\ \bosconfi{j}\\ i+1\notin\bosconfi{j}\\ i\in\bosconfi{i}}} c_{\bosconfi{i}\bosconfi{j}} a_{(\bosconfi{i}\cup\{i+1\})(\bosconfi{j}\cup\{i\})}.\label{rightsumbos}
\end{align}
Again, we observe that the set of monomials $\{a_{((\bosconfi{i}\setminus\{i\})\cup\{i+1\})\bosconfi{j}}\ |\ i+1\notin\bosconfi{j},\ i\in\bosconfi{i}\}$ is linearly independent for fixed $i$.
\par
By induction on the number $k_i\p$ of particles at position $i$ in $\bosconfi{i}$ we see that all coefficients  $c_{\bosconfi{i}\bosconfi{j}}$ are zero for $k_i\p\geq1$:\\
For $k_i\p=1$, the set $(\bosconfi{i}\setminus\{i\})\cup\{i+1\}$ does not contain any particle at position $i$ any more. Hence the monomial $a_{((\bosconfi{i}\setminus\{i\})\cup\{i+1\})\bosconfi{j}}$ cannot appear in the right sum in equation~\eqref{rightsumbos}, and its coefficient $c_{\bosconfi{i}\bosconfi{j}}$ must have been zero.
For the induction step we assume that the coefficients for all $a_{\bosconfi{i}\bosconfi{j}}$ with at most $k_i\p$ particles at position $i$ in the output particle configuration $\bosconfi{i}$ are zero. Consider some $a_{\bosconfi{i}\bosconfi{j}}$ with $k_i\p+1$ particles at position $i$ in $\bosconfi{i}$.
So the set $(\bosconfi{i}\setminus\{i\})\cup\{i+1\}$ contains $k_i\p$ particles at position $i$ in $\bosconfi{j}$, and the monomial $a_{((\bosconfi{i}\setminus\{i\})\cup\{i+1\})\bosconfi{j}}$ cannot appear in the sum~\eqref{leftsumbos}. Again we see that its coefficient $c_{\bosconfi{i}\bosconfi{j}}$ must have been zero.
\par
We have deduced now that only those monomials labelled by minimal input particle configurations $\bosconfi{j}=(r,0,\ldots,0)$ and output particle configuration $\bosconfi{i}=(0,\ldots,0,s)$ may have nonzero coefficients. Since the number of particles has to be the same in $\bosconfi{i}$ and $\bosconfi{j}$, any central element is of the form
$$\sum\limits_{r\in\ZZ_{\geq0}} c_{(0,\ldots,0,r)(r,0,\ldots,0)} a_{(0,\ldots,0,r)(r,0,\ldots,0)}$$
as claimed.
\end{bew}
\begin{bem}
In the proof of Theorem~\ref{theo:centpart} one has to be careful: One cannot simply compare the coefficients in equalities of the form 
$$a_i\left(\sum c_{\bosconfi{i}\bosconfi{j}}a_{\bosconfi{i}\bosconfi{j}}\right)=\left(\sum c_{\bosconfi{i}\bosconfi{j}}a_{\bosconfi{i}\bosconfi{j}}\right)a_i$$ since the partic algebra $\partic{N}$ has zero divisors, see Example \ref{ex:nottfree}.
Therefore, when we consider the coefficients $c_{\bosconfi{i}\bosconfi{j}}$, we first have to determine linearly independent sets of monomials, e.g. of the form
$$\{a_{((\bosconfi{i}\setminus\{i\})\cup\{i+1\})\bosconfi{j}}\ |\ i+1\notin\bosconfi{j},\ i\in\bosconfi{i}\}.$$
This is in fact an application of the faithfulness result from Theorem~\ref{theo:particfaith} combined with the normal form for monomials from Theorem~\ref{theo:particbasis}.
\end{bem}
\begin{bem}
The partic algebra is not finitely generated over its center: The center is concentrated in degree $\ZZ_{\geq0}\cdot(1,\ldots 1)$ with respect to the $\ZZ^{N-1}$-grading from Remark \ref{bem:particgrading}. On the other hand one can see from the normal form in Proposition~\ref{prop:nofo} that all $\ZZ_{\geq0}^{N-1}$-graded components of the partic algebra are nontrivial, hence the partic algebra cannot be finitely generated over its degree $\ZZ_{\geq0}\cdot(1,\ldots 1)$ component. 
\end{bem}

\section{A short comparison with the affine case}\label{sec:aff}
In this section we give a brief outlook to the affine case that will be treated in the followup work \cite{affplac}. We refrain from giving any details here since the computations are substantially harder, and there is constant danger of mixing the two cases as the differences can be quite subtle.
\par
An affine version of the plactic algebra is obtained by a very similar construction, except that the indices of the generators are now read modulo $N$.
The (local) affine plactic algebra $\aplac{N}$ is given in \cite[Definition 5.4]{ks} by the unital associative $\ground$-algebra generated by $a_0,a_1,\ldots,a_{N-1}$ subject to the affine plactic relations
\begin{eqnarray*}
 a_i a_j\ &=\ a_j a_i\quad &\text{for }i-j\neq\pm1\mod N,\\
 a_i a_{i-1} a_i\ &=\ a_i a_i a_{i-1}\quad &\text{for } i, i-1\in\ZZ/N\ZZ,\\
 a_i a_{i+1} a_i\ &=\ a_{i+1} a_i a_i\quad &\text{for } i, i+1\in\ZZ/N\ZZ.
\end{eqnarray*}
The affine plactic algebra acts on the polynomial ring $\ground[x_1,\ldots,x_N,q]$ in $N+1$ variables as follows:
\begin{align*}
 a_i\cdot x_1^{k_1}\ldots x_N^{k_N}q^t\ &=\ \begin{cases} x_1^{k_1}\ldots x_i^{k_i-1}x_{i+1}^{k_{i+1}+1}\ldots x_N^{k_N}q^t\quad &\text{if }k_i>0,\\ 0&\text{else,}\end{cases}\\
a_0\cdot x_1^{k_1}\ldots x_N^{k_N}q^t\ &=\ \begin{cases} x_1^{k_1+1}\ldots  x_N^{k_N-1}q^{t+1}\quad &\text{if }k_N>0,\\ 0&\text{else.}\end{cases}
\end{align*}
This representation is called the affine bosonic particle representation of the affine plactic algebra $\aplac{N}$.
Similar to the bosonic particle configurations in the classical case one can identify a monomial $x_1^{k_1}\ldots x_N^{k_N}$ with a particle configuration on a circle with $N$ positions, with $k_i$ particles lying at position $i$. The indeterminate $q$ protocols how often we apply $a_0$ to a particle configuration.
\begin{figure}[H]
\begin{center}
\begin{tikzpicture}[scale=0.6]

\begin{scope}[xshift=-7.5cm]
\path[use as bounding box] (-3.5,-1.5) rectangle (3.5,2.5);
\draw[thick] (0,0) circle (2cm);

\foreach \x in {45,90,...,360}
\draw[very thick] (\x:1.9cm) -- (\x:2.1cm);

\foreach \x in {0,...,7}
\draw ({90-\x*45}:1.5cm) node {$\x$};

\foreach \x in {0,1,2,5}
\fill ({90-\x*45}:2.35cm) circle (1.2mm);
\foreach \x in {1,5}
\fill ({90-\x*45}:2.7cm) circle (1.2mm);
\foreach \x in {1}
\fill ({90-\x*45}:3.05cm) circle (1.2mm);

\end{scope}

\begin{scope}[xshift=0cm]
\path[use as bounding box] (-3.5,-1.5) rectangle (3.5,2.5);
\draw[thick] (0,0) circle (2cm);

\foreach \x in {45,90,...,360}
\draw[very thick] (\x:1.9cm) -- (\x:2.1cm);

\foreach \x in {0,...,7}
\draw ({90-\x*45}:1.5cm) node {$\x$};

\foreach \x in {0,1,4,6}
\fill ({90-\x*45}:2.35cm) circle (1.2mm);
\foreach \x in {1,7}
\fill ({90-\x*45}:2.7cm) circle (1.2mm);
\foreach \x in {1}
\fill ({90-\x*45}:3.05cm) circle (1.2mm);

\foreach \x in {2,5}
\draw[thick, densely dotted] ({90-\x*45}:2.35cm) circle (1.2mm);
\foreach \x in {5}
\draw[thick, densely dotted] ({90-\x*45}:2.7cm) circle (1.2mm);

\foreach \x in {2}
\draw[thick, ->] ({85-\x*45}:2.35cm) arc ({85-\x*45}:{45-\x*45}:2.35cm);
\foreach \x in {3}
\draw[thick, ->] ({90-\x*45}:2.35cm) arc ({90-\x*45}:{50-\x*45}:2.35cm);

\foreach \x in {5}
\draw[thick, ->] ({85-\x*45}:2.35cm) arc ({85-\x*45}:{50-\x*45}:2.35cm);

\foreach \x in {5}
\draw[thick, ->] ({85-\x*45}:2.7cm) arc ({85-\x*45}:{45-\x*45}:2.7cm);
\foreach \x in {6}
\draw[thick, ->] ({90-\x*45}:2.7cm) arc ({90-\x*45}:{50-\x*45}:2.7cm);

\end{scope}

\begin{scope}[xshift=7.5cm]
\path[use as bounding box] (-3.5,-1.5) rectangle (3.5,2.5);
\draw[thick] (0,0) circle (2cm);

\foreach \x in {45,90,...,360}
\draw[very thick] (\x:1.9cm) -- (\x:2.1cm);

\foreach \x in {0,...,7}
\draw ({90-\x*45}:1.5cm) node {$\x$};

\foreach \x in {0,1,4,6,7}
\fill ({90-\x*45}:2.35cm) circle (1.2mm);
\foreach \x in {1}
\fill ({90-\x*45}:2.7cm) circle (1.2mm);
\foreach \x in {1}
\fill ({90-\x*45}:3.05cm) circle (1.2mm);

\end{scope}
\end{tikzpicture}
\end{center}
\caption{Example for $N=8$: Application of $a_6 a_5 a_3a_2a_5$ to the particle configuration $(3,1,0,0,2,0,0,1)$ gives $(3,0,0,1,0,1,1,1)$.}
\label{fig:bosonicaff}
\end{figure}
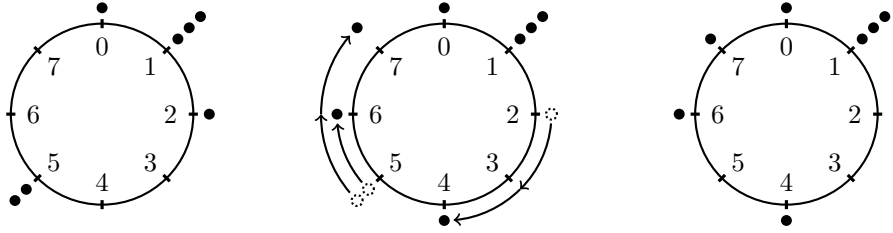
Again like in the classical case this action is not faithful. It factors over a quotient of $\aplac{N}$ by the relations
$$a_i a_{i-1} a_{i+1} a_i \ =\ a_{i+1} a_i a_{i-1} a_i\quad \text{for } i\in\ZZ/N\ZZ, $$
which is the affine version of the defining relation \eqref{part} of the partic algebra in the classical case. But it does not generate all the relations in the affine case: In addition there are infinitely many relations of the form
\begin{align*}
a_{i-1}^{m\p}a_i^m a_{i+1}^{k_{i+1}}a_{i+2}^{k_{i+2}}\ldots a_{i-2}^{k_{i-2}}a_{i-1}^m\ &=\ a_i^m a_{i-1}^{m\p}a_{i+1}^{k_{i+1}}a_{i+2}^{k_{i+2}}\ldots a_{i-2}^{k_{i-2}}a_{i-1}^m\quad &\text{for } i\in\ZZ/N\ZZ\\
a_i^m a_{i+1}^{k_{i+1}}a_{i+2}^{k_{i+2}}\ldots a_{i-2}^{k_{i-2}}a_{i-1}^m a_i^{m\p}\ &=\ a_i^m a_{i+1}^{k_{i+1}}a_{i+2}^{k_{i+2}}\ldots a_{i-2}^{k_{i-2}}a_i^{m\p}a_{i-1}^m \quad &\text{for } i\in\ZZ/N\ZZ,
\end{align*}
where $m, m\p, k_{i+1},k_{i+2},\ldots,k_{i-2}\in\ZZ_{\geq0}$ are nonnegative integers. These relations can be seen as the proper affine version of the partic relations \eqref{plac1}, \eqref{plac2} and \eqref{part}, since we find in \cite{affplac} that they generate the kernel of the action of $\aplac{N}$ on the affine bosonic particle representation $\ground[x_1,\ldots,x_N,q]$.


\begin{thebibliography}{XXXXX}

\bibitem[BFZ96]{bfz}
	{A. Berenstein, S. Fomin and A. Zelevinsky},
	{\it Parametrizations of canonical bases and totally positive matrices},
	Adv. Math. {\bf 122} (1996), no. 1, 49--149.

\bibitem[BJS93]{bjs}
	{S. C. Billey, W. Jockusch and R. P. Stanley},
	{\it Some combinatorial properties of {S}chubert polynomials},
	J. Algebraic Combin. {\bf 2} (1993), no. 4, 345--374.

\bibitem[BM16]{bejo}
  {G. Benkart and J. Meinel},
  {\it The center of the affine nil{T}emperley-{L}ieb algebra},
  Math. Z. {\bf 284} (2016), no. 1-2, 413--439.  

\bibitem[FG98]{fomingreene}
	{S. Fomin and C. Greene},
	{\it Noncommutative {S}chur functions and their applications},
	Discrete Math. {\bf 193} (1998), no. 1-3, 179--200.

\bibitem[Fom95]{fomin-schur}
  {S. Fomin},
	{\it Schur operators and {K}nuth correspondences},
	J. Combin. Theory Ser. A {\bf 72} (1995), no. 2, 277--292.

\bibitem[Ful97]{fulton-yt}
	{W. Fulton},
	{\it Young tableaux},
	London Mathematical Society Student Texts, vol. 35, Cambridge University Press, Cambridge, 1997.

\bibitem[Hey98]{hey}
	{Heyworth, A.},
	{\it Rewriting as a special case of non-commutative {G}r\"obner basis theory},
	Computational and geometric aspects of modern algebra ({E}dinburgh, 1998), London Math. Soc. Lecture Note Ser., vol. 275, Cambridge Univ. Press, Cambridge, 2000, pp. 101--105.

\bibitem[HK02]{hongkang}
	{J. Hong and S.-J. Kang},
	{\it Introduction to quantum groups and crystal bases},
	Graduate Studies in Mathematics, vol. 42, American Mathematical Society, Providence, RI, 2002.

\bibitem[KS10]{ks}
	{C. Korff and C. Stroppel},
	{\it The {$\h{\mathfrak{sl}}(n)_k$}-{WZNW} fusion ring: a combinatorial construction and a realisation as quotient of quantum cohomology},
	Adv. Math. {\bf 225} (2010), no. 1, 200--268.
	
\bibitem[LS81]{ls-plac}
	{A. Lascoux and M.-P. Sch{\"u}tzenberger},
	{\it Le monoide plaxique},
	Noncommutative structures in algebra and geometric combinatorics ({N}aples, 1978), Quad. ``Ricerca Sci.'', vol. 109, CNR, Rome, 1981, pp. 129--156.

\bibitem[Mei16]{diss}
	{J. Meinel},
	{\it Affine nil{T}emperley-{L}ieb Algebras and Generalized {W}eyl Algebras: {C}ombinatorics and Representation Theory},
	Dissertation, University of Bonn, 2016.
	
\bibitem[Mei]{affplac}
	{J. Meinel},
	{\it A plactic algebra action on bosonic particle configurations: {T}he affine case},
	In preparation.

\bibitem[Rei01]{reinge}
  {M. Reineke},
  {\it Generic extensions and multiplicative bases of quantum groups at {$q=0$}},
  Represent. Theory {\bf 5} (2001), 147--163 (electronic).
  
\bibitem[Rei02]{reinmono}
  {M. Reineke},
  {\it The quantic monoid and degenerate quantized enveloping algebras},
  arXiv/math/0206095 (2002).
  
\bibitem[Rin90]{ringel}
	{C. M. Ringel},
	{\it Hall algebras and quantum groups},
	Invent. Math. {\bf 101} (1990), no. 3., 583--591.
	    
\bibitem[Ste03]{stembridge}
	{J. R. Stembridge},
	{\it A local characterization of simply-laced crystals},
  Trans. Amer. Math. Soc. {\bf 355} (2003), no. 12, 4807--4823 (electronic).



\end{thebibliography}
\end{document}